\documentclass[titlepage,11pt]{article}
\usepackage{latexsym}
\usepackage{amsfonts}
\usepackage{amsmath}
\newcommand\blackslug{\hbox{\hskip 1pt \vrule width 4pt height 8pt depth 1.5pt
        \hskip 1pt}}
\newcommand\bbox{\hfill \quad \blackslug \bigbreak}

\def\l{,\ldots,}

\title{Rooted grid minors}
\author{
D\'aniel Marx\\
Computer and Automation Research Institute,\\
Hungarian Academy of Sciences (MTA SZTAKI)\\
\and 
Paul Seymour\thanks{Supported by ONR grant N00014-10-1-0680 and NSF grant DMS-0901075.}\\
Princeton University, Princeton, NJ 08544\\
\and 
Paul Wollan\thanks{Supported by the European Research Council under the European Unions Seventh Framework
Programme (FP7/2007-2013)/ERC Grant Agreement no. 279558.}\\
Department of Computer Science,\\ 
University of Rome, ``La Sapienza'', Rome, Italy.\\
}

\date{February 17, 2013; revised \today}

\newtheorem{thm}{}[section]

\newcommand{\Proof}{\noindent{\bf Proof.}\ \ }

\begin{document}
\maketitle
\begin{abstract}
Intuitively, a {\em tangle} of large order in a graph is a highly-connected part of the graph, and it is known that
if a graph has a tangle of large order then it has a large grid minor. Here we show that
for any $k$, if $G$ has a tangle of large order and $Z$ is a set of vertices of cardinality $k$ that cannot be separated from 
the tangle by any separation of order less than $k$, then $G$ has a large grid minor containing $Z$, in which the members of $Z$
all belong to the outside of the grid. This is a lemma for use in a later paper.
\end{abstract}

\section{Introduction}
A {\em separation} of {\em order $k$} in a graph $G$ is a pair $(A,B)$ of subgraphs of $G$ such that $A\cup B = G$, 
$E(A\cap B) = \emptyset$, and $|V(A\cap B)| = k$. A {\em tangle} in $G$ of {\em order $\theta \ge 1$} is a set $\mathcal{T}$
of separations of $G$, all of order less than $\theta$, such that
\begin{itemize}
\item for every separation $(A,B)$ of order less than $\theta$, $\mathcal{T}$ contains one of $(A,B),(B,A)$
\item if $(A_i,B_i)\in \mathcal{T}$ for $i = 1,2,3$, then $A_1\cup A_2\cup A_3\ne G$
\item if $(A,B)\in \mathcal{T}$ then $V(A) \ne V(G)$.
\end{itemize}
Let $G,H$ be graphs. A {\em pseudomodel} of $H$ in $G$ is a map $\eta$ with domain $V(H)\cup E(H)$, where
\begin{itemize}
\item for every $v\in V(H)$, $\eta(v)$ is a non-null subgraph of $G$, all pairwise vertex-disjoint
\item for every edge $e$ of $H$, $\eta(e)$ is an edge of $G$, all distinct
\item if $e\in E(H)$ and $v\in V(H)$ then $e\notin E(\eta(v))$
\item for every edge $e=uv$ of $H$, if $u\ne v$ then $\eta(e)$ has
one end in $V(\eta(u))$ and the other
in $V(\eta(v))$; and if $u = v$, then $\eta(e)$ is an edge of $G$ with all ends in $V(\eta(v))$.
\end{itemize}
If in addition we have
\begin{itemize}
\item $\eta(v)$ is connected for each $v\in V(H)$
\end{itemize} 
then we call $\eta$ a {\em model} of $H$ in $G$.
Thus, $G$ contains $H$ as a minor if and only if there is a model of $H$ in $G$.
If $\eta$ is a pseudomodel of $H$ in $G$, and $F\subseteq V(H)$, we denote 
$$\bigcup(V(\eta(v))\;:v\in F)$$
by $\eta(F)$; and if $F$ is a subgraph of $H$, $\eta(F)$ denotes the subgraph of $G$ formed by the union of all the subgraphs
$\eta(v)$ for $v\in V(F)$ and all the edges $\eta(e)$ for $e\in E(F)$.

For $g\ge 1$, the {\em $g\times g$-grid} has vertex set $\{v_{ij}\::1\le i,j\le g\}$, and vertices
$v_{ij}, v_{i'j'}$ are adjacent if $|i'-i|+|j'-j| = 1$. We denote this graph by $\mathcal{G}_g$.
For $1\le i\le g$, we call 
$\{v_{i1},v_{i2}\l v_{ig}\}$ a {\em row} of the grid, and define the {\em columns} of the grid similarly.

The following was proved in~\cite{GM5,GM10}:
\begin{thm}\label{gridthm}
For all $g\ge 1$ there exists $K\ge 1$ with the following property. Let $\mathcal{T}$ be a tangle of order at least $K$
in a graph $G$. Then there is a model $\eta$ of $\mathcal{G}_g$ in $G$, such that 
for each $(A,B)\in \mathcal{T}$, if $\eta(R)\subseteq V(A)$ for some row $R$ of the grid, then $(A,B)$ has order at least $g$.
\end{thm}

Our objective here is an analogous result, for graphs with some vertices distinguished, the following:

\begin{thm}\label{main}
For all $k,g$ with $1\le k\le g$ there exists $K\ge 1$ with the following property. Let $\mathcal{T}$ be a tangle of order at least $K$
in a graph $G$, and let $Z\subseteq V(G)$ with $|Z| = k$. Suppose that there is no separation $(A,B)\in \mathcal{T}$
of order less than $k$ with $Z\subseteq V(A)$. Then there is a model $\eta$ of $\mathcal{G}_g$ in $G$, such that 
\begin{itemize}
\item for $1\le i\le k$, $V(\eta(v_{i,1}))$ contains a member of $Z$
\item for each $(A,B)\in \mathcal{T}$, if $\eta(R)\subseteq V(A)$ for some row $R$ of the grid, then $(A,B)$ has order at least $g$.
\end{itemize}
\end{thm}

A form of this result is implicit in a paper of Bruce Reed (statement 5.5 of ~\cite{reed}), but what we need is not explicitly 
proved there, so it seems necessary to do it again.
It has as an immediate corollary the following (the proof of which is clear):

\begin{thm}\label{planarout}
Let $H$ be a planar graph, drawn in the plane, and let $v_1\l v_k$ be distinct vertices of $H$, each incident with the infinite region.
Then there exists $K$ with the following property. Let $\mathcal{T}$ be a tangle of order at least $K$ in a graph $G$, 
and let $Z\subseteq V(G)$ with $|Z| = k$ such that there is no separation $(A,B)\in \mathcal{T}$
of order less than $k$ with $Z\subseteq V(A)$. Then there is a model $\eta$ of $H$ in $G$ such that for $1\le i\le k$, $\eta(v_i)$
contains a vertex of $Z$.
\end{thm}

\section{The main proof}

To prove \ref{main}, it is convenient to prove something a little stronger, which we explain next. Let $H$ be a subgraph of $G$.
We define $\beta_G(H)$ to be the set of vertices of $H$ incident with an edge of $G$ that does not belong to $E(H)$, and call
$\beta_G(H)$ the {\em boundary} of $H$ in $G$. If $f\in E(G)$, $G/f$ denotes the graph obtained from $G$ by contracting $f$.

Let $G$ be a graph and $Z\subseteq V(G)$ with $|Z| = k$. Let $\eta$ be a model of $\mathcal{G}_g$ in $G$. 
We say $\eta$ is {\em $Z$-augmentable} in $G$ if there is a model $\eta'$ of $\mathcal{G}_g$ in $G$, and we can
label the vertices of $\mathcal{G}_g$ as usual, such that
\begin{itemize}
\item for $1\le i\le g$ and $2\le j\le g$, $\eta'(v_{ij}) = \eta(v_{ij})$
\item for $1\le i\le g$, $\eta'(v_{i1}) = \eta(v_{i1})$ if $i>k$, and $\eta'(v_{i1}) \supseteq \eta(v_{i1})$ if $i\le k$
\item for $1\le i\le k$, $V(\eta'(v_{i1}))$ contains a member of $Z$
\item for each $e\in E(\mathcal{G}_g)$, $\eta'(e) = \eta(e)$.
\end{itemize}
In this case we call $\eta'$ a {\em $Z$-augmentation} of $\eta$ in $G$.

\begin{thm}\label{gridbits}
Let $g,k$ be integers with $g\ge k\ge 1$, and let $n$ be an integer such that $n> k(g+2k)$.
Let $G$ be a graph, and let
$Z\subseteq V(G)$ with $|Z| = k$. Let $J$ be a subgraph of $\mathcal{G}_n$, with boundary $\beta$, including at least one row of $\mathcal{G}_n$.
Let $\eta$ be a pseudomodel of $J$ in $G$.
Suppose that
\begin{itemize}
\item[(i)] for each $v\in V(J)$, either $\eta(v)$ is connected and $v\notin \beta$, 
or every component of $\eta(v)$ contains a vertex of $Z$
\item[(ii)] there is no separation $(A,B)$ of $G$ of order less than $k$ such that $Z\subseteq V(A)$ and there is a row $R$ of $\mathcal{G}_n$
with $R\subseteq V(J)$ and $\eta(R)\subseteq V(B)$.
\end{itemize}
Then there is a subgraph $H$ of $J$, isomorphic to $\mathcal{G}_g$,
such that $Z\cap V(\eta(v))$ is null for every $v\in V(H)$, and the restriction of $\eta$ to $H$ is $Z$-augmentable.
\end{thm}
\Proof 
We proceed by induction on $|V(G)|+|E(G)|$. 
\\
\\
(1) {\em We may assume that there is no separation $(A,B)$ of $G$ of order $k$ with $B\ne G$ such that 
$Z\subseteq V(A)$ and there is a row $R$ of $\mathcal{G}_n$
with $R\subseteq V(H)$ and $\eta(R)\subseteq V(B)$.}
\\
\\
For suppose that $(A,B)$ is such a separation. Let $J'$ be the subgraph of $J$ with vertex set those $v\in V(J)$ with $\eta(v)\cap B$
non-null, and with edge set all edges $e$ of $J$ such that $\eta(u)\cap A$ is null for some end $u$ of $e$. (Note that
if $e\in E(J')$ then both ends of $e$ in $J$ belong to $V(B)$, so $J'$ is well-defined.) Let
$\beta'$ be the boundary of $J'$ in $\mathcal{G}_n$. By hypothesis $V(J')$ includes at least one row of $\mathcal{G}_n$. Let
$Z' = V(A\cap B)$. Then $|Z'| = k$. For each $v\in V(J')$, let $\eta'(v) = \eta(v)\cap B$ (note that $\eta(v)\cap B$ is non-null
from the definition of $V(J')$), and for each $e\in E(J')$,
let $\eta'(e) = \eta(e)$ (note that $\eta(e)\in E(B)$ from the definition of $E(J')$.) 
Thus $\eta'$ is a pseudomodel of $J'$ in $B$.

Now let $v\in V(J')$. We must show that either $\eta'(v)$ is connected and $v\notin \beta'$,
or every component of $\eta'(v)$ contains a vertex of $Z'$. 
We may assume that some component $C'$ of $\eta'(v)$ is disjoint from $Z'$. Let $C$ be the component of $\eta(v)$ containing $C'$.
If $C\ne C'$, then some vertex $u\in V(C')$ is adjacent in $C$ to some vertex $v\in V(C)\setminus V(C')$, and consequently 
$v\notin V(B)$; but then $u\in V(A\cap B) = Z'$, a contradiction. So $C=C'$. If some vertex of $C$ is in $Z$, then
that vertex belongs to $V(A)$ and hence to $Z'$, a contradiction. Thus no vertex of $C$ is in $Z$. It follows from hypothesis (i)
of the theorem that 
$\eta(v)$ is connected and $v\notin \beta$. In particular, since 
$$C'\subseteq \eta'(v)\subseteq \eta(v) = C = C'$$
it follows that $\eta'(v)$ is connected. It remains to check that
$v\notin \beta'$. Thus, suppose $v\in \beta'$, and so $v$ is incident in $\mathcal{G}_n$ with some edge $f = uv$ 
where $f\notin E(J')$.
Since $v\notin \beta$, it follows that $f\in E(J)$ and $u\in V(J)$. Since $f\notin E(J')$, both of $\eta(u),\eta(v)$ have
non-null intersection with $A$. But then $\eta(v)$ meets $Z'$, a contradiction. This proves that
for every $v\in V(J')$, either $\eta'(v)$ is connected and $v\notin \beta'$,
or every component of $\eta'(v)$ contains a vertex of $Z'$.

We claim that there is no separation $(A',B')$ of $B$ of order less than $k$ such that $Z'\subseteq V(A')$ 
and there is a row $R$ of $\mathcal{G}_n$
with $R\subseteq V(J')$ and $\eta(R)\subseteq V(B')$. For suppose there is such a separation $(A',B')$. 
Then $(A\cup A', B')$ is a separation of $G$. Moreover, $(A\cup A')\cap B' = A'\cap B'$, since
$$A\cap B'\subseteq A\cap B = Z'\subseteq A'.$$ 
But this contradicts hypothesis (ii) of the theorem.

Consequently, we may apply the inductive hypothesis with $G,Z,J,\eta, \beta$ replaced by $B,Z',J',\eta',\beta'$. We deduce that
there is a subgraph $H$ of $J'$, isomorphic to $\mathcal{G}_g$, such that $Z'\cap V(\eta'(v))$ is null for every $v\in V(H)$, 
and the restriction of $\eta'$ to $H$ is $Z$-augmentable in $B$.
Let $v\in V(H)$. Since $Z'\cap \eta'(v)$, it follows that $\eta'(v) = \eta(v)$. We deduce that
the restriction of $\eta$ to $H$ is $Z$-augmentable in $G$, and so the theorem holds. This proves (1).
\\
\\
(2) {\em We may assume that if $f\in E(G)$, then 
$f = \eta(e)$ for some $e\in E(J)$.} Consequently for each $v\in V(J)$, 
either $\eta(v)$ has only one vertex, or $V(\eta(v))\subseteq Z$.
\\
\\
For suppose not. Suppose first that there is no $u\in V(J)$ with $f\in E(\eta(u))$.
It follows that $\eta$ is a pseudomodel of $J$ in $G\setminus f$. By (1), hypothesis (ii) of the theorem holds for
$G\setminus f, Z,J,\eta, \beta$; and the other hypothesis holds trivially. Thus from the inductive hypothesis,
the theorem holds for $G\setminus f, Z,J,\eta, \beta$ and hence for $G, Z,J,\eta, \beta$.  
We may therefore assume that there exists $u\in V(J)$ with $f\in E(\eta(u))$. If $f$ is a loop or both ends of $f$ belong to $Z$, 
define $\eta'(u) = \eta(u)\setminus f$, and $\eta'(v) = \eta(v)$ for every other vertex $v$ of $J$; then $\eta'$ is
 a pseudomodel of $J$ in $G\setminus f$, and again the result follows from the inductive hypothesis. Finally, if $f$
is not a loop and some end of $f$ does not belong to $Z$, define
$\eta'(u) = \eta(u)/f$, and $\eta'(v) = \eta(v)$ for every other vertex $v$ of $J$; then $\eta'$ is
 a pseudomodel of $J$ in $G/f$, and again the result follows from (1) and the inductive hypothesis. This proves the
first assertion of (2), and the second follows.

\bigskip
Now let us label the vertices of $\mathcal{G}_n$ as usual. Let $Z'$ be the set of all vertices $v$ of $\mathcal{G}_n$ such that
$Z\cap \eta(v)\ne \emptyset$. Since $|Z|=k$ it follows that $|Z'|\le k$, and $\beta\subseteq Z'$ from hypothesis (i).
\\
\\
(3) {\em There is a subgraph $H_0$ of $J$, isomorphic to $\mathcal{G}_{g+2k}$, such that every row
of $\mathcal{G}_n$ that intersects $V(H_0)$ is a subset of $V(J)\setminus Z'$.}
\\
\\
From the choice of $n$, 
there are $k+1$ subgraphs of $\mathcal{G}_n$, each isomorphic to $\mathcal{G}_{g+2k}$, such that no row of 
$\mathcal{G}_n$ meets more than one of them. Consequently
there is a subgraph $H_0$ of
$\mathcal{G}_n$, isomorphic to $\mathcal{G}_{g+2k}$, such that no row of $\mathcal{G}_n$ meets both $V(H_0)$ and $Z'$. 
Let $H'$ be the subgraph of $\mathcal{G}_n$ induced on the union of the rows of $\mathcal{G}_n$ that meet $V(H_0)$.
We claim that every vertex of $H'$ belongs to $J$.
For suppose not; then none of them belong to $J$, since $H'$ is connected and none of its vertices belong to $\beta\subseteq Z'$.
Since there is a row $R$ of $\mathcal{G}_n$ with $R\subseteq V(J)$, it follows that every column of $\mathcal{G}_n$ meets both
$V(J)$ and $V(H')$,
and therefore meets $\beta$ and hence $Z'$. But $|Z'|\le k<n$, a contradiction. This proves (3).

\bigskip

For $1\le i\le n$, let $R_i = \{v_{ij}\;:1\le j\le n\}$; thus $R_i$ is a row of $\mathcal{G}_n$.
Since $V(H_0)\cap Z' = \emptyset$, (2) implies that $|V(\eta(v))| = 1$ for each $v\in V(H_0)$. Let 
$v'_{ij}$ be the unique vertex of $\eta(v_{ij})$ for all $v{ij}\in V(H_0)$.
Let the vertices of $H_0$ be the set of all $v_{ij}$ where $i_0-k\le i\le i_0+k+ g-1$ and $j_0-k\le j\le j_0+k+g-1$.  
Let $L$ be the set of all vertices $v_{ij_0}$ where $i_0\le i\le i_0+k-1$. Thus $|L| = k$. 
For $1\le s\le k$, let $H_s$ be the subgraph of $H_0$ induced on the vertex set
$$\{v_{ij}:s-k \le i-i_0,j-j_0\le g-1+k-s\},$$
and for $0\le s\le k-1$ let $C_i$ be the cycle of $H_0$ induced on the vertex set $V(H_s)\setminus V(H_{s+1})$.
Let $H = H_k$; thus
$H$ is isomorphic to $\mathcal{G}_g$.
Let $G^* = G\setminus \eta(V(H)\setminus L)$.
We may assume that
\\
\\
(4) {\em There is a separation $(A^*,B^*)$ of $G^*$ of order less than $k$, such that $Z\subseteq V(A^*)$ and $\eta(L)\subseteq V(B^*)$.}
\\
\\
For if not, then by Menger's theorem there are $k$ vertex-disjoint paths
$P_1\l P_k$ of $G^*$, where $v'_{i_0+i-1,j_0}$ belongs to $P_i$ for $1\le i\le k$. For each $v\in V(H)$, let
$\eta'(v) = \eta(v)$ if $v\notin L$, and $\eta'(v) = \eta(v)\cup P_i$ if $v \in L$ and $v = v_{i_0+i-1,j_0}$. Let $\eta(e) = \eta(e)$
for each edge $e$ of $H$. Then $\eta'$ is a $Z$-augmentation of $H$, and the theorem holds. This proves (4).
\\
\\
(5) {\em There exists $s$ with $0\le s\le k-1$ such that $\eta(C_s)\subseteq B^*$.}
\\
\\
For since $|A^*\cap B^*|<k$, there exists $s\in \{0\l k-1\}$ such that $\eta(C_s)\cap A^*\cap B^*$ is null, and hence
$\eta(C_s)$ is a subgraph of one of $A^*,B^*$. Suppose that $\eta(C_s)\subseteq A^*$. 
Now let $R$ be a row of $\mathcal{G}_n$ that meets $L$, and let $P$ be the path of $\mathcal{G}_n$
between $L$ and $V(C_s)$. Since $R\subseteq V(J)\setminus Z'$, and there

Suppose that $\eta(C_s)$ is a subgraph of $A^*$, and for $0\le i\le k-1$ let
$P_i$ be the path of $H_0$ between $L$ and $V(C_s)$ with vertex set included in $R_{i_0+i}$. It follows that some vertex
of $P_i$ is in $V(A^*)$ (namely its end in $\eta(L)$)
\bigskip

Let $X=V(A\cap B)$. Let $A',B',X'$ be respectively the sets of vertices $v$
of $J$ satisfying $\eta(v)\cap V(A)\ne \emptyset$, $\eta(v)\cap V(B)\ne \emptyset$, and $\eta(v)\cap X\ne \emptyset$. 
We claim that 
\\
\\
(5) {\em The following hold:
\begin{itemize}
\item If $v\in A'$, then every component of $\eta(v)$ contains a vertex of $V(A)$.
\item If $v\in B'\setminus A'$, then $\eta(v)$ is connected.
\item $A'\cap B' = X'$.
\item If $a'\in A'\setminus B'$ and $b'\in B'\setminus A'$ then $a',b'$ are not adjacent in 
$J$. 
\item If $C$ is a connected subgraph of $\mathcal{G}_n$ disjoint from $X'$ and with non-empty intersection with $B'$ then $C$ is a subgraph of $J$
and $V(C)\subseteq B'\setminus A'$. 
\end{itemize}
}
\noindent
For the first bullet, let $v\in A'$; the assertion is true if $\eta(v)$
is connected, and otherwise every component of $\eta(v)$ contains a vertex of $Z\subseteq V(A)$ as required.
For the second bullet, let $v\in B'\setminus A'$; then $Z\cap \eta(v) = \emptyset$, and so $\eta(v)$ is connected.
For the third bullet, clearly $X'\subseteq A'\cap B'$. For the converse, 
let $v\in A'\cap B'$, and choose $b\in V(B)\cap \eta(v)$. By the first bullet, the component
of $\eta(v)$ containing $b$ has a vertex in $V(A)$, and therefore a vertex in $X$,
since every path between $V(A), V(B)$ in $G^*$ contains a vertex of $X$; 
and therefore $v\in X'$. 

For the fourth bullet, suppose that $a'\in A'\setminus B'$ and $b'\in B'\setminus A'$ are adjacent in $J$, joined by an edge $f'$.
Let $\eta(f') = f$ say; then $f$ has an end in $\eta(a')$ and an end
in $\eta(b')$. Let $C$ be the component of $\eta(a')$ containing an end of $f$. By the first two bullets, the subgraph
formed by the union of $C$, $\eta(b')$, and $f$ is connected, and since it meets both $V(A)$ and $V(B)$, it also meets $X$,
and so one of $a',b'\in X'$, contrary to the third bullet.

Finally, for the fifth bullet, let $C$ be a connected subgraph of $\mathcal{G}_n$ disjoint from $X'$ and 
with non-empty intersection with $B'$. If the claim does not hold, then since 
$V(C)\cap A'\cap B' = \emptyset$ (by the third bullet), there are adjacent vertices $a',b'$ of $C$
with $b'\in B'\setminus A'$ and
$a'\in (A'\setminus B')\cup (V(\mathcal{G}_n)\setminus V(J))$. 
By the fourth bullet, $a'\notin A'$, and so $a'\notin V(J)$, and consequently
$b'\in \beta\subseteq Z'\subseteq A'$, and so $b'\in A'\cap B'$, a contradiction. This proves (5).

\bigskip

Since $|X'|<k$, there exists $r$ with $i_0\le r\le i_0+k-1$ such that $R_r\cap X' = \emptyset$. It follows from the fifth bullet of (5) 
that
$v_{rj}\in B'$ for $1\le j\le j_0$, since $L\subseteq B'$.

For $1\le s\le k$, let $S_s$ be the set of all $v_{i,j}$ where $(i,j)$ belongs to
\begin{eqnarray*}
     \{(i,j) :& i_0-k+s-1\le i\le i_0+k+g-s   ,& j\in \{j_0-k+s-1, j_0+k+g-s\}\} \\ 
\cup \{(i,j) :& i\in \{i_0-k+s-1, i_0+k+g-s\} ,& j_0-k+s-1\le j\le j_0+k+g-s\}.
\end{eqnarray*}
Thus, for $1\le s\le k$, $S_s$ is the vertex set of a cycle of $H_0$ ``surrounding'' $H$;  
and the sets $S_1\l S_k$ are pairwise disjoint
and each is disjoint from $V(H)$. 
Since $|X'|<k$, there exists $s$ with $1\le s\le k$ such that $S_s\cap X' = \emptyset$. 
Since $v_{r,j_0-k+s-1}\in B'\cap S_s$ it follows from the fifth bullet of (5) that $S_s\subseteq B'$.
\\
\\
(6) {\em There is a path of $G$ between $Z$ and $\eta(v_{r,j_0})$ disjoint from $X$.}
\\
\\
Suppose first that $R_r\cap Z'\ne \emptyset$, and let $P$ be a minimal subpath of $\mathcal{G}_n$ between $Z'$ and $v_{r,j_0}$ 
with $V(P)\subseteq R_r$. It follows that no vertex of $P$ except possibly one end belongs to $\beta$, since $\beta\subseteq Z'$;
and so $P$ is a path of $J$, and $\eta(v)$ is defined for every vertex $v$ of $P$, and therefore
the desired path can be chosen in $G|\eta(P)$. We may therefore assume that 
$R_r\cap Z' = \emptyset$, and so $R_r\subseteq V(J)$.
By hypothesis, there is no separation $(C,D)$ of $G$ of order less than $k$ such that 
$Z\subseteq V(C)$ and $V(D)$ includes $\eta(R_r)$. In particular there is a path $T$ of $G\setminus X$ between $Z$ and $\eta(R_r)$,
since $|X|<k$. But then the union of $T$ and $G|\eta(R_r)$ includes the required path. This proves (6).

\bigskip

Let $Y'$ be the union of $S_{s+1}\l S_k$ and $V(H)$; that is, the set of vertices of $\mathcal{G}_n$ ``surrounded'' by $S_s$.
By (6), there is a minimal path $Q$ of $G\setminus X$ between $Z$ and $\eta(Y')$; let its ends be
$z\in Z$ and $y\in \eta(Y')$. It follows that no vertex of $Q\setminus y$ is in $\eta(V(H)\setminus L)$, and hence $Q\setminus y$
is a path of $G^*$.
Let $y\in V(\eta(y'))$; then $y'\in Y'$. 
Let $x$ be the neighbour of $y$ in $Q$, and let $x\in \eta(x')$. From (2), the edge $xy$ of $G$ equals
$\eta(f')$ for some edge $f'$ of $J$ incident with $x',y'$, and since $x'\notin Y'$, it follows that $x'\in S_s$. Consequently
$Q\setminus y$ is a path of $G^*$ between $Z$ and $\eta(S_s)$ disjoint from $X$. Since
$(A,B)$ is a separation of $G^*$, and $\eta(S_s)\subseteq B$, it follows that $\eta(S_s)\cap X$ is non-null, a contradiction.
This proves \ref{gridbits}.~\bbox

Finally, let us deduce \ref{main}, which we restate:

\begin{thm}\label{main2}
For all $k,g$ with $1\le k\le g$ there exists $K\ge 1$ with the following property. Let $\mathcal{T}$ be a tangle of order at least $K$
in a graph $G$, and let $Z\subseteq V(G)$ with $|Z| = k$. Suppose that there is no separation $(A,B)\in \mathcal{T}$
of order less than $k$ with $Z\subseteq V(A)$. Then there is a model $\eta$ of $\mathcal{G}_g$ in $G$, such that
\begin{itemize}
\item for $1\le i\le k$, $V(\eta(v_{i,1}))$ contains a member of $Z$
\item for each $(A,B)\in \mathcal{T}$, if $\eta(R)\subseteq V(A)$ for some row $R$ of the grid, then $(A,B)$ has order at least $g$.
\end{itemize}
\end{thm}
\Proof
Let $n$ be as in \ref{gridbits}. Choose $K$ to satisfy \ref{gridthm} (with $g$ replaced by $n$.) We claim that this choice of $K$
satisfies \ref{main2}. For let $\mathcal{T}$ be a tangle of order at least $K$
in a graph $G$, and let $Z\subseteq V(G)$ with $|Z| = k$. Suppose that there is no separation $(A,B)\in \mathcal{T}$
of order less than $k$ with $Z\subseteq V(A)$. By \ref{gridthm} 
there is a model $\eta$ of $\mathcal{G}_n$ in $G$, such that
for each $(A,B)\in \mathcal{T}$, if $\eta(R)\subseteq V(A)$ for some row $R$ of $\mathcal{G}_n$, then $(A,B)$ has order at least $n$.
\\
\\
(1) {\em There is no separation $(A,B)$ of $G$ of order less than $k$ such that $Z\subseteq V(A)$ and there is a row $R$ of $\mathcal{G}_n$
with $R\subseteq V(J)$ and $\eta(R)\subseteq V(B)$.}
\\
\\
For suppose that $(A,B)$ is such a separation. Since $k\le n\le K$, it follows that one of $(A,B),(B,A)\in \mathcal{T}$.
But there is no separation $(A,B)\in \mathcal{T}$
of order less than $k$ with $Z\subseteq V(A)$, so $(A,B)\notin \mathcal{T}$;  
and for each $(C,D)\in \mathcal{T}$, if $\eta(R)\subseteq V(C)$ for some row $R$ of $\mathcal{G}_n$, then $(C,D)$ has order at least $n$,
so $(B,A)\notin \mathcal{T}$, a contradiction. This proves (1).

\bigskip
From (1) and \ref{gridbits}, taking $J = \mathcal{G}_n$, we deduce that
there is a subgraph $H$ of $\mathcal{G}_n$, isomorphic to $\mathcal{G}_g$,
such that the restriction of $\eta$ to $H$ is $Z$-augmentable. 
\\
\\
(2) {\em For each $(A,B)\in \mathcal{T}$, if $\eta(R)\subseteq V(A)$ for some row $R$ of $\mathcal{G}_g$, 
then $(A,B)$ has order at least $g$.}
\\
\\
For since $\eta(R)\subseteq V(A)$, and $J$ is a subgraph of $\mathcal{G}_n$, 
it follows that there are at least $g$ columns $C$ of $\mathcal{G}_n$ such that $C\cap V(A)\ne \emptyset$. If each of them
contains a vertex of $A\cap B$ then $|A\cap B|\ge g$ as required, and otherwise some column $C$ of $\mathcal{G}_n$
is included in $V(A)$. But then every row of $\mathcal{G}_n$ contains a vertex in $V(A)$; if they all meet
$A\cap B$ then $|A\cap B|\ge n\ge g$ as required, and otherwise some row of $\mathcal{G}_n$ is included in $V(A)$. But then
from the choice of $\eta$, $(A,B)$ has order at least $n\ge g$. This proves (2).

\bigskip

This proves \ref{main2}.~\bbox


\begin{thebibliography}{99}
\bibitem{reed} B.Reed, ``Mangoes and blueberries'', {\em Combinatorica} 19 (1999), 267--296.
\bibitem{GM5} N.Robertson and P.D.Seymour, ``Graph minors. V. Excluding a planar graph'',
{\em J. Combinatorial Theory, Ser. B,} 41 (1986), 92--114.
\bibitem{GM10} N.Robertson and P.D.Seymour, ``Graph minors. X. Obstructions to tree-decomposition'',
{\em J. Combinatorial Theory, Ser. B,} 52 (1991), 153--190.

\end{thebibliography}
\end{document}